\input amstex
\documentstyle{amsppt}
\TagsOnRight \NoRunningHeads
\magnification=\magstep1 \vsize 22 true cm \hsize 16 true cm

\centerline{\bf ON THE AUTOMORPHISMS OF SOME ONE-RELATOR GROUPS}

\bigskip

\centerline{\bf D.~Tieudjo and D.~I.~Moldavanskii}
\medskip

\eightpoint

\centerline{University of Ngaoundere}

\centerline{P.~o.~Box 454 Ngaoundere, Cameroon.}

\centerline{{\it E-mail:} tieudjo\@yahoo.com}
\medskip
\centerline{Ivanovo State University}

\centerline{Ermak str. 37, 153025 Ivanovo, Russia.}

\centerline{{\it E-mail:} moldav\@ivanovo.ac.ru}

\bigskip

{\noindent \bf Abstract:} The description of the automorphism
group of group $\langle a, b;\ [a^m,b^n]=1\rangle$ ($m,n>1$) in
terms of generators and defining relations is given. This result
is applied to prove that any normal automorphism of every such
group is inner.

\smallskip


\smallskip


\tenpoint

\bigskip\bigskip

\centerline{\bf Introduction}
\medskip

The automorphism group of certain one-relator groups was studied
by several authors. D.~Collins in [3] obtained the presentation by
generators and defining relations of the automorphism group of
Baumslag - Solitar groups $G(l,m)=\langle a, b;\
a^{-1}b^la=b^m\rangle$ when $|l|=1$ or $|m|=1$ or $|l|>1$, $|m|>1$
and $l$ and $m$ are coprime; in particular, in these cases the
group $\text{Aut}\,G(l,m)$ turns out to be finitely related. Later
D.~Collins and F.~Levin [4] found the presentation of the group
$\text{Aut}\,G(l,m)$ when $m=ls$, $|l|>1$ and $|s|>1$ and showed
thereby that in this case the group $\text{Aut}\,G(l,m)$ is not
finitely generated. In the same paper the more extensive class of
groups $G=\langle a_1, a_2, \dots, a_n, t;\ t^{-1}w^lt=w^m\rangle$
where $w$ is a word in $a_1$, $a_2$, \dots, $a_n$ was considered.
When $n\geqslant 2$, $w$ is neither a proper power nor primitive
in the free group $\langle a_1, a_2, \dots, a_n\rangle$ and $m=ls$
with $|s|>1$, authors gave the presentation of group
$\text{Aut}\,G$ and this group turns out to be infinitely
generated too. Some HNN-extension of Baumslag - Solitar groups
$G(l,m;k)=\langle a, t;\ t^{-1}a^{-k}ta^lt^{-1}a^kt=a^m\rangle$
were considered by A.M.Brunner in [2]. In the case when $|l|\neq
|m|$ he described all endomorphisms of such groups  and noted that
if $|l|=1$ or $|m|=1$ then the group $\text{Aut}\,G(m,n;k)$ is not
finitely generated. Using the Brunner's results, M.~Kavutskii and
D.~Moldavanskii [5] under assumption $|l|\neq |m|$ obtained the
presentation of $\text{Aut}\,G(m,n;k)$ and proved that this group
is finitely generated if and only if none of the integers $l$ and
$m$ is divisor of another. Furthermore, if the group
$\text{Aut}\,G(m,n;k)$ is finitely generated then it is finitely
related. It should be mentioned here that up to now it is unknown
whether the automorphism group of any one-relator group is
finitely presented if it is finitely generated.

Results listing above are relative to one-relator groups which in
either case are connected with Baumslag - Solitar groups. In the
present paper we consider another class of one-relator groups
consisting of groups $G_{mn}$ with presentation
$$
G_{mn}=\langle a, b;\ [a^m,b^n]=1\rangle
$$
where $m$ and $n$ are arbitrary integers satisfying inequalities
$m>1$ and $n>1$. We obtain the presentation of group
$\text{Aut}\,G_{mn}$ by generators and defining relations and
thereby prove that it is finitely related. We prove also that any
normal automorphism of every group $G_{mn}$ is inner.
\smallskip

As can be immediately verified the following mappings of
generators of group $G_{mn}$ define the automorphisms of $G_{mn}$
(which will be denoted by the same symbols):
$$
\align
&\lambda: a\mapsto a^{-1}, \quad b\mapsto b; \\
&\mu: a\mapsto a, \qquad b\mapsto b^{-1}; \\
&\nu: a\mapsto a^{-1}, \quad b\mapsto b^{-1}.\endalign
$$
It is evident that $\lambda^2=\mu^2=1$, $\lambda\mu=\mu\lambda$
and $\lambda\mu=\nu$ and therefore these automorphisms together
with the identity mapping constitute the subgroup $K$ of group
$\text{Aut}\,G_{mn}$ and $K$ is the Klein four-group. If $m=n$ the
mapping
$$
\eta: a\mapsto b, \quad b\mapsto a
$$
defines one more automorphism of $G_{mn}$. The relations
$\eta^2=1$, $\eta^{-1}\lambda\eta=\mu$ and
$\eta^{-1}\mu\eta=\lambda$ (which can also be immediately checked)
show that the subgroup $L$ of $\text{Aut}\,G_{mn}$ generated by
subgroup $K$ and element $\eta$ is the split extension of $K$ by
the 2-cycle $\langle\eta\rangle$. It will be shown here that if
$m\neq n$ then $\text{Aut}\,G_{mn}= K\cdot\text{Inn}\,G_{mn}$ and
if $m=n$ then $\text{Aut}\,G_{mn}= L\cdot\text{Inn}\,G_{mn}$. More
explicitly, we shall prove the following

\proclaim{Theorem 1} Let $\lambda$, $\mu$ and $\eta$ be the
automorphisms of group $G_{mn}$ defined above and $\alpha$ and
$\beta$ be the inner automorphisms of $G_{mn}$ generated by
elements $a$ and $b$ respectively.

If $m\neq n$ then group $\text{Aut}\,G_{mn}$ is generated by the
automorphisms $\lambda$, $\mu$, $\alpha$ and $\beta$ and defined
by the relations
$$
\aligned
&\text{1}.\ \lambda^2=\mu^2=1; \\
&\text{2}.\ \lambda\mu=\mu\lambda; \\
&\text{3}.\ \lambda^{-1}\alpha\lambda=\alpha^{-1}; \\
&\text{4}.\ \lambda^{-1}\beta\lambda=\beta; \endaligned \qquad
\qquad \qquad \qquad \qquad \qquad \aligned
&\text{5}.\ \mu^{-1}\alpha\mu=\alpha; \\
&\text{6}.\ \mu^{-1}\beta\mu=\beta^{-1}; \\
&\text{7}.\ \alpha^m\beta^n=\beta^n\alpha^m. \endaligned
$$

If $m=n$ then group $\text{Aut}\,G_{mn}$ is generated by the
automorphisms $\lambda$, $\mu$, $\eta$, $\alpha$ and $\beta$ and
defined by the relations 1 -- 7 and the additional relations
$$
\aligned
&\text{8}.\ \eta^2=1; \\
&\text{9}.\ \eta^{-1}\lambda\eta=\mu; \endaligned \qquad \qquad
\qquad \qquad \qquad \qquad \aligned &\text{10}.\
\eta^{-1}\alpha\eta=\beta. \endaligned
$$
\endproclaim

Theorem 1 can be applied to characterize the normal automorphisms
of groups $G_{mn}$. Let us recall that an automorphism of a group
$G$ is said to be normal if it maps onto itself every normal
subgroup of $G$. It is evident that any inner automorphism is
normal. In general, the converse is not true. It was proved in [6,
7] that any normal automorphism of a non-cyclic free group must be
inner. Generalizing this result M.~Neschadim [10] exhibited that
the same assertion is true for any group which is a non-trivial
free product. Also he gave the example of one-relator group
possessing a normal automorphism which is not inner. Nevertheless,
for groups $G_{mn}$ we have :

\proclaim{Theorem 2} Any normal automorphism of group $G_{mn}$ is
inner.
\endproclaim

We note that the residual finiteness of group $G_{mn}$ (i.~e.,
recall, for any non-identity element $g\in G_{mn}$ there exists a
homomorphism $\varphi$ of group $G_{mn}$ onto some finite group
$X$ such that $g\varphi\neq 1$) is well known; it follows, for
example, from the result of paper [1]. Since $G_{mn}$ is finitely
generated, then by Mal'cev theorem [9], it is Hopfian; i.e. every
of its surjective endomorphism is an automorphism. Some other
properties of these groups were considered in [12] where, in
particular, their construction as amalgamated free product and the
description of their endomorphisms were given. These results can
be used for somewhat shortening of the proof of our theorem 1 but
for completeness we shall give here the independent proof.
\bigskip

\centerline{\bf 1. Preliminaries}
\medskip

As we have just mentioned, the group $G_{mn}$ can be constructed
as amalgamated free product and we begin from some properties of
this group-theoretic construction.

Let $G=(A*B;\ H)$ be a free product of groups $A$ and $B$ with
amalgamated subgroup $H$. Then any element $g\in G$ can be written
in the form $g=x_1x_2 \cdots x_s$, where elements $x_1$, $x_2$,
\dots , $x_s$ belong in turns to one of groups $A$ and $B$ and if
$s>1$ then no one of them belongs to subgroup $H$. Such
representation is called a reduced form of element $g$ and the
number $s$ of factors of it (uniquely determined by $g$) is called
a length of $g$ and denoted by $l(g)$. An element $g$ is said to
be cyclically reduced if either $l(g)=1$ or the factors $x_1$ and
$x_s$ of its reduced form $g=x_1x_2 \cdots x_s$ do not belong to
the same subgroup $A$ or $B$ (the definition is correct since all
reduced forms of element $g$ have or do not have this property
simultaneously). It is easy to see that any element of $G$ is
conjugate with a cyclically reduced element. Moreover, an
immediate induction gives the

\proclaim{Proposition 1.1} If element $g$ of group $G=(A*B;\ H)$
is not cyclically reduced and $l(g)>1$ then $g$ can be written in
the form
$$
g=u\cdot v\cdot u^{-1},
$$
where elements $u$ and $v$  have reduced forms $u=x_1x_2 \cdots
x_r$ and $v=y_1y_2 \cdots y_s$ with $r\geqslant 1$ and $s\geqslant
1$, element $v$ is cyclically reduced, elements $x_r$ and $y_1$ do
not belong to the same subgroup $A$ or $B$ and if $s>1$ then
element $y_sx_r^{-1}$ does not belong to subgroup $H$.
\endproclaim

By means of proposition 1.1 it is easy to prove the

\proclaim{Proposition 1.2} Let element $g$ of group $G=(A*B;\ H)$
do not belong to subgroup $A$ and let $g^k\in A$ for some integer
$k\neq 0$. Then $g=x^{-1}yx$ for some $x,y\in G$ where element $y$
belongs to one of subgroups $A$ or $B$ and $y^k\in H$.
\endproclaim

Also we need the simple

\proclaim{Proposition 1.3} Let $G=(A*B;\ H)$ and suppose that the
amalgamated subgroup $H$ is contained in the centre of both groups
$A$ and $B$. If element $g\in G$ does not belong to subgroup $A$
then $g^{-1}Ag\cap A=H$.
\endproclaim

Indeed, the inclusion $H\subseteq g^{-1}Ag\cap A$ is evident. To
prove the inverse inclusion let $\rho$ be the natural homomorphism
of group $G$ onto quotient group $G/H$ which is the ordinary free
product of quotients $A/H$ and $B/H$. Then since $g\rho\notin
A\rho$ we have
$$
(g^{-1}Ag\cap A)\rho\subseteq (g\rho)^{-1}(A\rho)(g\rho)\cap
(A\rho)=1
$$
and therefore $g^{-1}Ag\cap A\subseteq H$.
\smallskip

Further, we need the construction of group $G_{mn}$ in terms of
amalgamated free product. For this purpose let $H=\langle c, d;\
[c,d]=1\rangle$ be the free abelian group of rank 2, $A=(\langle a
\rangle * H;\ a^m=c)$ be the amalgamated free product of infinite
cycle $\langle a \rangle$ and $H$ and $B=(H * \langle b \rangle;\
d=b^n)$ be the amalgamated free product of $H$ and infinite cycle
$\langle b \rangle$. Then it easy to show by means of Tietze
transformations that group $G_{mn}$ is isomorphic to the free
product $(A*B;\ H)$ of groups $A$ and $B$ with amalgamated
subgroup $H$. These notations are assumed in what follows.

Since in constructions of groups $A$ and $B$ the amalgamated
subgroups are central in the free factors, proposition 1.3 gives
the

\proclaim{Proposition 1.4} If an element $g$ of group $A$ does not
belong to subgroup $H$ then $g^{-1}Hg\cap H=\langle c\rangle$, and
if an element $g$ of group $B$ does not belong to subgroup $H$
then $g^{-1}Hg\cap H=\langle d\rangle$.
\endproclaim

\proclaim{Proposition 1.5} Any element $g$ of group $G_{mn}$ such
that $g^{-1}Hg\cap H\neq 1$ is contained in subgroup $A$ or in
subgroup $B$.
\endproclaim

For the proof it is enough to show that if $g=x_1x_2 \cdots x_s$
is reduced form of $g$ with $s>1$ then $g^{-1}Hg\cap H = 1$. Let
us suppose that $x_1\in A$; the case $x_1\in B$ is considered
similarly. For any element $h\in H$ the inclusion $g^{-1}hg\in H$
implies the inclusions $x_1^{-1}hx_1\in H$ and
$x_2^{-1}(x_1^{-1}hx_1)x_2\in H$. Thus, $x_1^{-1}hx_1\in
x_1^{-1}Hx_1\cap H$ and since $x_1\in A\setminus H$ it follows
from proposition 1.4 that $x_1^{-1}hx_1=c^k$ for some integer $k$.
Similarly, inclusion $x_2^{-1}c^kx_2\in x_2^{-1}Hx_2\cap H$ gives
$x_2^{-1}c^kx_2=d^l$ for some integer $l$, and since element $d^l$
lies in the centre of group $B$, we have the equality $c^k=d^l$.
As elements $c$ and $d$ form the basis of free abelian group $H$,
hence $k=l=0$ and $h=1$.
\smallskip

\proclaim{Proposition 1.6} Any abelian subgroup of group $G_{mn}$
which contains a cyclically reduced element of length greater than
1 is cyclic.
\endproclaim

\demo{Proof} Let $U$ be abelian subgroup of group $G_{mn}$ and let
$U$ contain a cyclically reduced element $u$ of length greater
than 1. It is not difficult to see that any element of $G_{mn}$
commuting with $u$ is either element of $H$ or cyclically reduced
of length greater then 1. Since proposition 1.5 implies $U\cap
H=1$ we conclude that all nonidentity elements of $U$ are
cyclically reduced of length greater than 1.

Let $u$ be the nonidentity element of $U$ of the smallest length
and $u=u_1u_1 \cdots u_r$ be a reduced form of it. We claim that
subgroup $U$ is generated by $u$. Namely, for any nonidentity
element $v\in U$ we shall prove by induction on $l(v)$ that $v$ is
equal to some power of $u$.

Let $v=v_1v_2 \cdots v_s$ be a reduced form of element $v$.
Replacing, if necessary, element $v$ by element $v^{-1}$, we can
assume that elements $u_1$ and $v_1$ belong to the same subgroup
$A$ or $B$. Then since elements $v_s$ and $u_1$ do not belong to
the same subgroup $A$ or $B$ and the right side of equation
$$
u_r^{-1} \cdots u_2^{-1}u_1^{-1}v_1v_2 \cdots v_s u_1u_2 \cdots
u_r =v_1v_2 \cdots v_s
$$
is cyclically reduced the product $h=u_r^{-1} \cdots
u_2^{-1}u_1^{-1}v_1v_2 \cdots v_r$ must be element of $H$. So, if
$s=r$ we have $h=u^{-1}v\in U$ and since $U\cap H=1$ we obtain the
equality $v=u$ giving the basis of induction.

If $s>r$ then $v=uv^{\prime}$, where $v^{\prime}=hv_{r+1} \cdots
v_s$. Since $l(v^{\prime})< s$ then by induction $v^{\prime}=u^k$
for some integer $k$. Hence $v=u^{k+1}$ and proof is complete.
\enddemo
\bigskip

\centerline{\bf 2. Proof of Theorem 1}
\medskip

\proclaim{Proposition 2.1} For any automorphism $\varphi$ of group
$G_{mn}$ there exists an inner automorphism $\psi$ of $G_{mn}$
such that either $a(\varphi\psi)\in A$ and $b(\varphi\psi)\in B$
or $a(\varphi\psi)\in B$ and $b(\varphi\psi)\in A$.
\endproclaim

\demo{Proof} Let $\varphi$ be an automorphism of group $G_{mn}$
and $u=a\varphi$, $v=b\varphi$. At first, we note that elements
$u$ and $v$ cannot be cyclically reduced of length greater than 1.

If, on the contrary, element $u$ is cyclically reduced and
$l(u)>1$ then element $u^m$ is also cyclically reduced of length
greater than 1, and since $[u^m,v^n]=1$, by proposition 1.6
elements $u^m$ and $v^n$ generate the (infinite) cyclic subgroup.
Therefore, $u^{mr}=v^{ns}$ for some nonzero integers $r$ and $s$.
But this equation implies the equation $a^{mr}=b^{ns}$ which is
not satisfied in group $G_{mn}$.

On the other hand, element $u$ is conjugate with a cyclically
reduced element and after multiplying $\varphi$ by suitable inner
automorphism we can assume that $u$ is cyclically reduced.
Consequently, by the remark above $u\in A$ or $u\in B$.

Suppose firstly that $u\in A$.  We claim that if $v\notin B$ then
$v=xyx^{-1}$ where $x\in A$ and $y\in B$ and therefore
$x^{-1}ux\in A$ and $x^{-1}vx\in B$. So, multiplying $\varphi$ by
one more inner automorphism we obtain the desired result.

Since $u$ and $v$ generate the group $G_{mn}$, $v\notin A$. Hence
if $v\notin B$ then $l(v)>1$ and since $v$ is not cyclically
reduced it has by proposition 1.1 the form
$$
v=x_1x_2 \cdots x_r\cdot y_1y_2 \cdots y_s\cdot (x_1x_2 \cdots
x_r)^{-1}
$$
where $r\geqslant 1$, $s\geqslant 1$, element $x_1x_2 \cdots x_r$
is reduced, element $y_1y_2 \cdots y_s$ is cyclically reduced,
elements $x_r$ and $y_1$ do not belong to the same subgroup $A$ or
$B$ and if $s>1$ then element $y_sx_r^{-1}$ does not belong to
subgroup $H$.

We assert now that the assumption $x_1\in B$ leads to the
contradiction. To prove this, let us note firstly that if $x_1\in
B$ then $l(v^n)>1$ and the first syllable of reduced form of $v^n$
is $x_1$. This is evident if $s>1$ or if $s=1$ and $y_1^n\notin
H$. If $s=1$ then $y_1\in B$ since if $y_1\in A$ then elements $u$
and $v$ are contained in the normal closure in $G_{mn}$ of
subgroup $A$ and therefore cannot generate the group $G_{mn}$.
Hence $x_r\in A$ and $r>1$. If $y_1^n\in H$ then $y_1^n\in
y_1^{-1}Hy_1\cap H$ and by proposition 1.4 $y_1^n=d^k$ for some
integer $k\neq 0$. Therefore $x_ry_1^nx_r^{-1}\in A\setminus H$,
$l(v^n)=2r-1>1$ and the first syllable of reduced form of $v^n$ is
$x_1$.

Now, since $x_1\in B$ the equality $v^{-n}u^mv^n=u^m$ implies
inclusions $u^m\in H$ and $x_1^{-1}u^mx_1\in H$. Since $u\in
A\setminus H$ (because the quotient group of $G_{mn}$ by the
normal closure of $H$ is not cyclic) and $x_1\in B\setminus H$ the
proposition 1.4 implies that $u^m=1$, a contradiction.

So, $x_1\in A$. If $v$ does not have the form claimed above then
$r>1$ and elements $u_1=x_1^{-1}ux_1$ and $v_1=x_1^{-1}vx_1$ turn
out to be in the previous case.

Thus, we have proved that if element $u=a\varphi$ belongs to
subgroup $A$ then after multiplying, if necessary, automorphism
$\varphi$ by one more inner automorphism we have $u\in A$ and
$v\in B$. Similar arguments will show that if element $u=a\varphi$
belongs to subgroup $B$ then after multiplying, if necessary,
automorphism $\varphi$ by one more inner automorphism we get $u\in
B$ and $v\in A$.
\enddemo

\proclaim{Proposition 2.2} Let elements $u$ and $v$ of group
$G_{mn}$ be such that $u\in A\setminus H$, $v\in B\setminus H$ and
$[u^r,v^s]=1$ for some integers $r\neq 0$ and $s\neq 0$. Then
$u=x^{-1}a^kx$ and $v=y^{-1}b^ly$ where $x\in A$, $y\in B$ and
nonzero integers $k$ and $l$ are such that $kr$ is divided by $m$
and $ls$ is divided by $n$.
\endproclaim

\demo{Proof} We note, firstly, that $u^r\in H$ and $v^s\in H$.
Indeed, if, say, $u^r\notin H$ then since $u^r\in A$, $v^s\in B$
and $[u^r,v^s]=1$ we get $v^s\in H$. Hence $v^s=u^{-r}v^su^r\in
H\cap u^{-r}Hu^r$ and $v^s=v^{-1}v^sv\in H\cap v^{-1}Hv$.
Proposition 1.4 implies now that $v^s=1$ which is impossible.

Since $u\in A\setminus H$ and $u^r\in H$ then proposition 1.2,
applied to the group $A=(\langle a \rangle * H;\ a^m=c)$, gives
$u=x^{-1}zx$ where $x\in A$, element $z$ is contained in subgroup
$\langle a \rangle$ or in subgroup $H$ and element $z^r$ belongs
to subgroup $\langle a^m \rangle$. But if $z\in H$ then the
inclusion $z^r\in \langle c \rangle$ is possible only if $z\in
\langle c \rangle$. Thus, in any case $z=a^k$ for some integer
$k\neq 0$ and $m$ divides $kr$ because $z^r\in\langle a^m
\rangle$.

So, we have proved that $u$ has the required form. The assertion
on the element $v$ is proved similarly.
\enddemo

\proclaim{Proposition 2.3} Let $F$ be the subgroup of group
$G_{mn}$ generated by elements $x^{-1}a^kx$ and $y^{-1}b^ly$ where
$x\in A$, $y\in B$ and $k,l\in \Bbb Z$. Then $F=G_{mn}$ if and
only if $|k|=1=|l|$ and $x\in \langle a\rangle \cdot\langle
d\rangle$, $y\in \langle b\rangle \cdot\langle c\rangle$.
\endproclaim

\demo{Proof} If $|k|=1=|l|$ and $x=a^pd^q$, $y=b^rc^s$ for some
integers $p,\ q,\ r,\ s$ then elements $c=(x^{-1}ax)^m$ and
$d=(y^{-1}by)^n$ belong to subgroup $F$ and since
$a=d^q(x^{-1}ax)d^{-q}$ and $b=c^s(y^{-1}by)c^{-s}$ we have $a\in
F$ and $b\in F$ and therefore $F=G_{mn}$.

Conversely, let us suppose that $F=G_{mn}$. Then the quotient
group of $G_{mn}$ by its commutator subgroup $G'_{mn}$ is
generated by elements $a^k$ and $b^l$ and since $G_{mn}/G'_{mn}$
is the free abelian group with basis $a$, $b$ we must have
$|k|=1=|l|$.

Let $A_1$ denote the subgroup of $G_{mn}$ generated by subgroup
$H$ and element $x^{-1}ax$ and let $B_1$ denote the subgroup of
$G_{mn}$ generated by subgroup $H$ and element $y^{-1}by$. Since
$H\leqslant A_1\leqslant A$ and $H\leqslant B_1\leqslant B$, it
follows by the theorem of H.~Neumann (see e.~g.~ [11], p.~512)
that the subgroup $F_1$ generated by $A_1$ and $B_1$ is the free
product of groups $A_1$ and $B_1$ with amalgamated subgroup $H$
and $A\cap F_1=A_1$, $B\cap F_1=B_1$. Therefore, since $F\leqslant
F_1$ the equality $F=G_{mn}$ implies $A_1=A$ and $B_1=B$.

Now we shall prove that if $A_1=A$ then $x\in \langle a\rangle
\cdot\langle d\rangle$. Let $x=x_1x_2\cdots x_r$ be the reduced
form of element $x$ (in the decomposition of group $A$ in
amalgamated product $A=(\langle a \rangle * H;\ a^m=c)$).

If $x\notin \langle a\rangle \cdot\langle d\rangle$ then
$r\geqslant 2$ and if $r=2$ then $x_1\in H$ and $x_2\in \langle
a\rangle$. If $r>2$ and $x_1\in \langle a\rangle$ then letting
$x'=x_2x_3\cdots x_r$ we see that subgroup $A_1$ is generated by
subgroup $H$ and element $(x')^{-1}a(x')$. Now if $x_r\in H$ then
$r>3$ and letting $x''=x_2x_3\cdots x_{r-1}$ we see that subgroup
$A_1$ is generated by subgroup $H$ and element $(x'')^{-1}a(x'')$.
Thus, we have shown that if $x\notin \langle a\rangle \cdot\langle
d\rangle$ then we can assume without loss of generality that
$r\geqslant 2$ and $x_1\in H$, $x_r\in \langle a\rangle$.

Let $\overline A$ be the quotient of group $A$ by central subgroup
$\langle a^m\rangle$ and $\overline g$ denote the image of element
$g\in A$ under the natural homomorphism of $A$ to $\overline A$.
Then $\overline A$ is the ordinary free product of cyclic group
$\langle \overline a\rangle$ of order $m$ and of infinite cycle
$\langle \overline d\rangle$. The image $\overline A_1$ of
subgroup $A_1$ is generated by element $\overline d$ and by image
$\overline{x^{-1}ax}$ of element $x^{-1}ax$. Our assumptions about
$x$ imply that
$$
\overline x_r^{-1}\cdots \overline x_2^{-1}\overline x_1^{-1}
\overline a \, \overline x_1\overline x_2\cdots \overline x_r
$$
is the reduced form of element $\overline{x^{-1}ax}$. This in turn
implies that any alternating product of nonidentity powers of
elements $\overline{x^{-1}ax}$ and $\overline d$ is reduced as
written. Thus, $\overline A_1\neq \overline A$ (since $\overline
a\notin \overline A_1$) and hence $A_1\neq A$. Consequently, the
equality $A_1=A$ really implies the inclusion $x\in \langle
a\rangle \cdot\langle d\rangle$ and the same arguments will show
that the equality $B_1=B$ implies the inclusion $y\in \langle
b\rangle \cdot\langle c\rangle$.
\enddemo

Now we can complete the proof of Theorem 1. Let $\varphi$ be an
automorphism of group $G_{mn}$. Proposition 2.1 implies that for
some inner automorphism $\psi$ of $G_{mn}$ we shall get either
$a(\varphi\psi)\in A$ and $b(\varphi\psi)\in B$ or
$a(\varphi\psi)\in B$ and $b(\varphi\psi)\in A$.

Firstly, let us consider the case when $a(\varphi\psi)\in A$ and
$b(\varphi\psi)\in B$. Since elements $a(\varphi\psi)$ and
$b(\varphi\psi)$ generate the group $G_{mn}$ and hence no one of
them belong to subgroup $H$, it follows from proposition 2.2 that
$a(\varphi\psi)=x^{-1}a^kx$ and $b(\varphi\psi)=y^{-1}b^ly$ for
some $x\in A$, $y\in B$ and nonzero integers $k$ and $l$. Now,
proposition 2.3 implies that
$a(\varphi\psi)=d^{-p}a^{\varepsilon}d^p$ and
$b(\varphi\psi)=c^{-q}b^{\delta}c^q$ for some integers $p$ and $q$
and $\varepsilon , \delta =\pm 1$. Then the product of
$\varphi\psi$ by the inner automorphism generated by element
$c^{-q}d^{-q}$ belongs to subgroup $K$ and therefore $\varphi\in
K\cdot\text{Inn}\,G_{mn}$.

Now, let $a(\varphi\psi)\in B$ and $b(\varphi\psi)\in A$. Then by
proposition 2.2 $a(\varphi\psi)=y^{-1}b^ly$ and
$b(\varphi\psi)=x^{-1}a^kx$ for some $x\in A$, $y\in B$ and
nonzero integers $k$ and $l$ where $kn$ is divided by $m$ and $lm$
is divided by $n$. Since proposition 2.3 again gives $|k|=1=|l|$,
conditions of divisibility imply the equality $m=n$. Thus, if
$m\neq n$ then $\text{Aut}\,G_{mn}= K\cdot\text{Inn}\,G_{mn}$.

If $m=n$ then the group $G_{mn}$ has the automorphism $\eta$ and
since $A\eta =B$ and $B\eta =A$ we obtain $a(\varphi\psi\eta)\in
A$ and $b(\varphi\psi\eta)\in B$. Therefore, automorphism
$\varphi\psi\eta$ belongs to subgroup $K\cdot\text{Inn}\,G_{mn}$.
This means that $\varphi\in L\cdot\text{Inn}\,G_{mn}$. Thus, in
the case $m=n$ we obtain $\text{Aut}\,G_{mn}=
L\cdot\text{Inn}\,G_{mn}$.

The validity of relations 1 -- 10 in the statement of theorem 1
can be checked immediately (and this in part was singled out
above) and it remains to show that these relations do define the
group $\text{Aut}\,G_{mn}$. Making use of relations 3 -- 6 in the
case $m\neq n$ and of relations 3 -- 6 and 10 in the case $m=n$,
any relation in the pointed out generators of $\text{Aut}\,G_{mn}$
can be transformed to the form $uv=1$ where $u$ is a product of
elements $\lambda$ and $\mu$ (or $\lambda$, $\mu$ and $\eta$) and
$v$ is a product of elements $\alpha$ and $\beta$. Since the unit
is the only element of subgroups $K$ and $L$ inducing the identity
automorphism of quotient group $G_{mn}/G'_{mn}$, we can conclude
that
$$
K\cap \text{Inn}\,G_{mn}=1 \qquad \text{and} \qquad L\cap
\text{Inn}\,G_{mn}=1
$$
and therefore the relation $uv=1$ implies $u=1$ and $v=1$. Since
relations 1 and 2 define the group $K$ and relations 1, 2, 8 and 9
define the group $L$, the relation $u=1$ is derivable from the
relations singled out in Theorem. Since the presentation above of
group $G_{mn}$ as amalgamated free product with regard to
corollary 4.5 in [8] makes evident the triviality of its centre,
the group $\text{Inn}\,G_{mn}$ is isomorphic to $G_{mn}$ and
therefore the relation $v=1$ must be derivable from the relation
7. Thus, any relation in the indicated generators of group
$\text{Aut}\,G_{mn}$ is derivable from the relations 1 -- 10 and
the proof is complete.
\bigskip

\centerline{\bf 3. Proof of Theorem 2}
\medskip

We begin with the rather obvious remark. If $\varphi$ is a normal
automorphism of a group $G$ and if $N$ is a normal subgroup of
group $G$ then the mapping $\overline\varphi$ of the factor group
$G/N$ onto itself, defined by
$$
(gN) \overline\varphi = (g\varphi) N \qquad (g\in G),
$$
is an automorphism of group $G/N$ and this automorphism is normal
too. The automorphism $\overline\varphi$ is said to be induced by
automorphism $\varphi$.

Now, let $\varphi$ be a normal automorphism of group $G_{mn}$.
Then by Theorem 1 $\varphi=\xi\psi$ where
$\psi\in\text{Inn}\,G_{mn}$ and $\xi\in K$ if $m\neq n$ or $\xi\in
L$ if $m=n$. Since automorphism $\varphi$ is normal if and only if
the automorphism $\xi$ is normal, it remains to show that any
non-identity element of subgroups $K$ and $L$ is not normal
automorphism.

Let $M$ and $N$ denote the normal closure in group $G_{mn}$ of
elements $a^m$ and $b^n$ respectively. Then the quotient group
$G_{mn}/M$ is the free product of cycle $\langle a \rangle$ of
order $m$ and infinite cycle $\langle b \rangle$ and the quotient
group $G_{mn}/N$ is the free product of infinite cycle $\langle a
\rangle$ and cycle $\langle b \rangle$ of order $n$.

Since the orders of elements $aM$ and $bM$ of the group $G_{mn}/M$
are different, then any automorphism of form $\kappa\eta$ where
$\kappa\in K$ does not induce any automorphism of this quotient
and therefore is not normal by the remark above.

In the same quotient group $G_{mn}/M$ the elements $bM$ and
$(bM)^{-1}$ are not conjugate, since two elements of a free factor
of an ordinary free product are conjugate if and only if they are
conjugate in the factor. Therefore, automorphisms $\overline\mu$
and $\overline\nu$ of group $G_{mn}/M$, induced by the
automorphisms $\mu$ and $\nu$ respectively, are not inner and
consequently, by the mentioned above result in [10],
$\overline\mu$ and $\overline\nu$ are not normal. Hence, from the
remark above it follows that automorphisms $\mu$ and $\nu$ of
group $G_{mn}$ are not normal. Analogously, automorphism $\lambda$
induces a non-inner automorphism in the quotient $G_{mn}/N$ and
therefore is not normal. Theorem 2 is demonstrated.
\bigskip

\centerline{\bf References}
\medskip

\item{[1]} G.~Baumslag, Free subgroups of certain one-relator
groups defined by positive words, {\it Math. Proc. Camb. Phil.
Soc.} {\bf 93} (1983) 247-251.

\item{[2]} A.~M.~Brunner, On a class of one-relator groups, {\it
Can. J. Math.} {\bf 32} 2 (1980) 414-420.

\item{[3]} D.~J.~Collins, The automorphism towers of some
one-relator groups, {\it Proc. London Math. Soc.} {\bf 36} (1978)
480-493.

\item{[4]} D.~J.~Collins and F. Levin, Automorphisms and hopficity
of certain Baumslag-Solitar groups, {\it Arch. Math.} {\bf 40}
(1983) 385-400.

\item{[5]} M.~Kavutskii and D.~Moldavanskii, On the class of
one-relator groups, {\it in  "Algebraic and discrete systems".
Ivan. State Univ.} (1988) 35-48. (Russian)

\item{[6]} A.~Lubotski, Normal automorphisms of free groups, {\it
J. Algebra} {\bf 63} 2 (1980) 494-498.

\item{[7]} A.~S.-T.~Lue, Normal automorphisms of free groups {\it
J. Algebra} {\bf 64} 1 (1980) 52-53.

\item{[8]} W.~Magnus, A.~Karrass and D.~Solitar, {\it
Combinatorial group theory}, John Wiley and Sons, Inc. New York
London Sydney, 1966.

\item{[9]} A.~I.~Mal'cev, On isomorphic representation of infinite
groups by matrixes {\it Math. Sbornik} {\bf 8} (1940) 405--422.
(Russian)

\item{[10]} M.~V.~Neschadim, Free products of groups do not have
outer normal automorphisms,  {\it Algebra and Logic} {\bf 35} 5
(1996) 562-566, (Russian).

\item{[11]} B~H. Neumann, An assay on free products of groups with
amalgamations, {\it Phil. Trans. Royal Soc. of London} {\bf 246}
(1954) 503--554.

\item{[12]} D. Tieudjo and D.~I. Moldavanskii, Endomorphisms of
the group $G_{mn}=\langle a, b;$ \linebreak $[a^m,b^n]=1\rangle$
($m,n>1$), {\it Afrika Matematika, J. of African Math. Union}
Series 3 {\bf 9} (1998) 11--18.

\end